\documentclass[12pt]{article}
\usepackage{latexsym}
\usepackage{amsfonts}
\usepackage{amssymb}
\usepackage{amscd}
\usepackage{array}
\usepackage{amsmath}
\usepackage{epsfig}
\usepackage{a4}
\begin{document}
\title{\textbf{Stringy $E$-functions of hypersurfaces and of Brieskorn singularities}\\ \date{}}
\author{Jan Schepers and Willem Veys\footnote{Partially supported by project G.0318.06 from the Research Foundation - Flanders (FWO). During the preparation of the manuscript the first-named author was partially supported by the Research Foundation - Flanders (FWO) and partially by the Netherlands Organisation for Scientific Research (NWO).}}
\maketitle
\begin{center}
\footnotesize{\textbf{Abstract}}
\end{center}
{\footnotesize We show that for a hypersurface Batyrev's stringy $E$-function can be seen as a residue of the Hodge zeta function, a specialization of the motivic zeta function of Denef and Loeser. This is a nice application of inversion of adjunction. If an affine hypersurface is given by a polynomial that is non-degenerate with respect to its Newton polyhedron, then the motivic zeta function and thus the stringy $E$-function can be computed from this Newton polyhedron (by work of Artal, Cassou-Nogu\`es, Luengo and Melle based on an algorithm of Denef and Hoornaert). We use this procedure to obtain an easy way to compute the contribution of a Brieskorn singularity to the stringy $E$-function. As a corollary, we prove that stringy Hodge numbers of varieties with a certain class of strictly canonical Brieskorn singularities are nonnegative. We conclude by computing an interesting 6-dimensional example. It shows that a result, implying nonnegativity of stringy Hodge numbers in lower dimensional cases, obtained in our previous paper, is not true in higher dimension.}\\

\section{Introduction} 
\noindent \textbf{1.1.} In \cite{Batyrev}, Batyrev defined the stringy $E$-function, an interesting singularity invariant of complex algebraic varieties with at most log terminal singularities. He used this function to formulate a topological mirror symmetry test for Calabi-Yau varieties with singularities, thereby extending the classical mirror symmetry test for smooth Calabi-Yau's. Let us recall Batyrev's definition, the related notion of `stringy Hodge numbers' and a remarkable conjecture stated by Batyrev.\\

\noindent \textbf{1.2.} The Grothendieck group of complex algebraic varieties, denoted
$K_0(Var_{\mathbb{C}})$, is the abelian group generated by the
symbols $[X]$, where $X$ is a complex algebraic variety (not
necessarily irreducible), and with the following relations:
\begin{itemize}
\item if $X$ is isomorphic to $Y$, then $[X]=[Y]$,
\item if $Y$ is a Zariski closed subset of $X$, then $[X]=[X\setminus Y] + [Y]$.
\end{itemize}
There is a product structure making $K_0(Var_{\mathbb{C}})$ into a
ring, defined by $[X]\cdot [Y]=[X\times Y]$. Thus the Grothendieck
ring is the value ring of the `universal Euler characteristic' on
algebraic varieties. The class of the affine line $\mathbb{A}^1_{\mathbb{C}}$
is usually denoted by $\mathbb{L}$; the class of a point is the
unity 1. \\

\noindent For an arbitrary complex variety $X$ of dimension $d$, the compactly supported cohomology $H_c^{\bullet}(X)$ carries a natural mixed Hodge structure, see \cite{Deligne1} and \cite{Deligne2} (we always use cohomology with complex coefficients). The data of this mixed Hodge structure are encoded in the Hodge-Deligne polynomial, defined by
\[ H(X;u,v):= \sum_{p,q=0}^d \left[\sum_{i=0}^{2d} (-1)^i h^{p,q}(H^i_c(X))\right] u^pv^q,\] where $h^{p,q}(H^i_c(X))$ denotes the dimension of the $(p,q)$-component $H^{p,q}(H_c^i(X))$. It is well known that the Hodge-Deligne polynomial induces a ring morphism from $K_0(Var_{\mathbb{C}})$ to $\mathbb{Z}[u,v]$. It maps $\mathbb{L}$ to $uv$. Note that $H(X;1,1)=\chi_{\text{top}}(X)$.\\

\noindent \textbf{1.3.} A normal irreducible variety $Y$ is called $\mathbb{Q}$-Gorenstein if a multiple $rK_Y$ of its canonical divisor is Cartier for some $r\in \mathbb{Z}_{>0}$ (we call $Y$ Gorenstein if $K_Y$ itself is Cartier). For example, all hypersurfaces and more generally all complete intersections are Gorenstein. Let $Y$ be $\mathbb{Q}$-Gorenstein and $f:X\to Y$ be a log resolution (this is a proper birational morphism from a smooth variety $X$ such that the exceptional locus of $f$ is a divisor with smooth irreducible components $D_i,i\in I,$ and normal crossings). Then $rK_X-f^*(rK_Y) = \sum_i b_iD_i$, with the $b_i\in \mathbb{Z}$. This is also formally written as $K_X - f^*(K_Y) = \sum_i a_iD_i$, with $a_i=b_i/r$. The variety $Y$ is called terminal, canonical, log terminal or log canonical if $a_i > 0, a_i \geq 0,a_i > -1$ or $a_i\geq -1$, respectively, for all $i$ (this does not depend on the chosen log resolution). We say that $Y$ is strictly canonical if it is canonical but not terminal. The number $a_i$ is called the discrepancy coefficient of $D_i$ and the difference $K_X - f^*(K_Y)$ is called the discrepancy. These definitions play a key r\^{o}le in the Minimal Model Program.\\

\noindent \textbf{1.4.} Now we are ready to define Batyrev's stringy $E$-function (see \cite{Batyrev}). Let $Y$ be an irreducible complex variety with at most log terminal singularities. Take a log resolution $f:X\to Y$ and denote the irreducible components of the exceptional locus by $D_i,i\in I$. For a subset $J\subset I$ write $D_J:=\cap_{i\in J} D_i$ and $D_J^{\circ}:=D_J\setminus \cup_{i\in I\setminus J} D_i$ ($D_{\emptyset}$ is taken to be $X$).
The \textit{stringy $E$-function} of $Y$ is
\[E_{st}(Y;u,v):= \sum_{J\subset I} H(D_J^{\circ};u,v) \prod_{i\in J}
\frac{uv-1}{(uv)^{a_i+1}-1}, \] where $a_i$ is the discrepancy
coefficient of $D_i$ and where the product $\prod_{i\in J}$ is 1
if $J=\emptyset$. Batyrev proved that this definition is independent of the chosen
log resolution (\cite[Theorem 3.4]{Batyrev}). His proof uses motivic integration. Alternatively, one can use the Weak Factorization Theorem by Abramovich, Karu, Matsuki and W{\l}odarczyk from \cite{weakfactorization}.\\

\noindent \textbf{1.5. Remark.} \begin{itemize}
\item[(1)] The stringy $E$-function can obviously be written as $H(Y_{ns};u,v)\,+$\,con\-tri\-bu\-tions of the singularities, where $Y_{ns}$ denotes the nonsingular part of $Y$.
\item[(2)] If $Y$ is smooth, then $E_{st}(Y)=H(Y)$ and if $Y$
admits a crepant resolution $f:X \to Y$ (i.e.\ such that the discrepancy is 0), then
$E_{st}(Y)=H(X)$.
\item[(3)] If $Y$ is Gorenstein, then all
$a_i\in \mathbb{Z}_{\geq 0}$ and $E_{st}(Y)$ becomes a rational
function in $u$ and $v$. It is then an element of
$\mathbb{Z}[[u,v]]\cap \mathbb{Q}(u,v)$.
\item[(4)] The \textit{stringy Euler number} of $Y$ is defined as
\[\lim_{u,v \to 1} E_{st}(Y;u,v)=\sum_{J\subset I}
\chi_{\text{top}}(D_J^{\circ})\prod_{j\in J} \frac{1}{a_j+1}.\]
\item[(5)] It is easy to deduce the following alternative expression for the stringy $E$-function:
\[E_{st}(Y) = \sum_{J\subset I} H(D_J;u,v ) \prod_{i\in J} \frac{uv-(uv)^{a_i+1}}{(uv)^{a_i+1}-1}.\]
\end{itemize}

\noindent \textbf{1.6.} Assume moreover that $Y$ is projective of
dimension $d$. Then Batyrev proved the following instance of
Poincar\'e and Serre duality (\cite[Theorem 3.7]{Batyrev}):
\begin{itemize}
\item[(i)] $E_{st}(Y;u,v)=(uv)^dE_{st}(Y;u^{-1},v^{-1})$ (note that the Hodge-Deligne polynomial of a \textit{smooth} projective variety satisfies the same relation),
\item[(ii)]$E_{st}(Y;0,0)=1$.
\end{itemize}
If in addition $Y$ has at worst Gorenstein canonical singularities and if
$E_{st}(Y;u,v)$ is a polynomial $\sum_{p,q} a_{p,q}u^pv^q$, he
defined the \textit{stringy Hodge numbers} of $Y$ as
$h_{st}^{p,q}(Y):=(-1)^{p+q}a_{p,q}$. It is easy to see that
\begin{itemize}
\item[(1)] they can only be nonzero for $0\leq p\leq d$ and $0\leq q\leq d$,
\item[(2)] $h_{st}^{0,0}(Y)=h_{st}^{d,d}(Y)=1$,
\item[(3)] $h_{st}^{p,q}(Y)=h_{st}^{q,p}(Y)=h_{st}^{d-p,d-q}(Y)=h_{st}^{d-q,d-p}(Y)$,
\item[(4)] if $Y$ is smooth, the stringy Hodge numbers are equal
to the usual Hodge numbers.
\end{itemize}
\noindent So the stringy Hodge numbers satisfy many of the properties of the usual Hodge numbers of smooth projective varieties. There is though one desired property that is not clear at all: nonnegavity\,!\\

\noindent \textbf{Conjecture} \cite[Conjecture 3.10]{Batyrev}\textbf{.} \textsl{Stringy Hodge numbers are nonnegative.}\\

\noindent The idea is that stringy Hodge numbers should be dimensions of certain (pieces of) cohomology spaces, just like classical Hodge numbers. In specific cases conjectural definitions of `string cohomology' were given by Borisov and Mavlyutov in \cite{BorisovMavlyutov}. The authors also made connections with the orbifold cohomology of Chen and Ruan from \cite{ChenRuan}. The above conjecture is trivially true for all varieties that admit a crepant resolution (more specifically this is the case for all Gorenstein canonical surfaces). In \cite{SchepersVeys} the conjecture was proved for threefolds, and for varieties with isolated singularities that admit a log resolution with all discrepancy coefficients $>\lfloor \frac{d-4}{2} \rfloor$, where $d$ is the dimension. (So for $d=4,5$ this just means isolated terminal singularities.) In fact we proved a stronger statement. One can look at the power series development $\sum_{i,j\geq 0} b_{i,j}u^iv^j$ of the stringy $E$-function (so we do not assume that it is a polynomial), and for the mentioned cases we proved that $(-1)^{i+j}b_{i,j} \geq 0$ for $i+j\leq d$. In view of property (3) of stringy Hodge numbers above, this implies the conjecture if the stringy $E$-function is a polynomial.\\

\noindent \textbf{1.7.} In this paper we first show that the stringy $E$-function of a hypersurface can be computed as a kind of residue of the Hodge zeta function, which is a specialization of the motivic zeta function of Denef and Loeser (Section 2). In fact this is a nice application of a version of inversion of adjunction by Stevens (see \cite{Stevens}). This can be used to compute the stringy $E$-function for non-degenerate hypersurfaces since Artal, Cassou-Nogu\`es, Luengo and Melle describe in \cite{Spanjaarden} an algorithm to compute the motivic zeta function in this case (based on an algorithm by Denef and Hoornaert for Igusa's $p$-adic zeta function from \cite{DenefHoornaert}), see Section 3. We use this method to describe an easy way to compute the contribution of a Brieskorn singularity to the stringy $E$-function (Theorem 4.2). As a corollary, we prove the nonnegativity of the stringy Hodge numbers for varieties with certain strictly canonical Brieskorn singularities (Corollary 4.4). We conclude by computing an interesting example of a 6-dimensional variety (Section 5). It shows that the result of \cite{SchepersVeys} mentioned at the end of 1.6 is no longer true for isolated terminal singularities in dimension 6 (in this example the coefficient of $(uv)^3$ in the power series development of the stringy $E$-function is negative).

\section{The motivic zeta function, the Hodge zeta function and the stringy $E$-function}
\noindent \textbf{2.1.} Let $X$ be a smooth irreducible complex algebraic variety and let $f:X\to \mathbb{A}^1_{\mathbb{C}}$ be a non-constant morphism. In this context Denef and Loeser define the \textit{naive motivic zeta function} of $f$ (see \cite[Definition 3.2.1]{DenefLoeserBarcelona}; we just call it the motivic zeta function). For our purposes we do not need the original definition in terms of jet spaces of $X$, but only the formula in terms of an embedded resolution of $f^{-1}(0)$. Let $h: Z\to X$ be an embedded resolution of $f^{-1}(0)$. So $Z$ is a nonsingular
variety, $h$ is a proper birational morphism, the restriction $h:
Z\setminus h^{-1}(f^{-1}(0))\to X\setminus f^{-1}(0)$ is an
isomorphism and $h^{-1}(f^{-1}(0))$ is a divisor with smooth
irreducible components and normal crossings on $Z$. Denote the
irreducible components of $h^{-1}(f^{-1}(0))$ by $E_i, i\in I$,
and for a subset $J\subset I$, use the notations $E_J $
and $E_J^{\circ}$ as in the introduction. Let $N_i$ be the multiplicity of $E_i$ in the divisor
of $f\circ h$ and let $\nu_i-1$ be the multiplicity of $E_i$ in
the divisor of $h^{*}dx$, where $dx$ is a local generator of the
sheaf of differential forms of maximal degree on $X$. These two
numbers are called the \textit{numerical data} of $E_i$. Denote the localization of the Grothendieck ring $K_0(Var_{\mathbb{C}})$ with respect to $\mathbb{L}$ by $\mathcal{M}_{\mathbb{C}}$. The motivic zeta function is then the following element of $\mathcal{M}_{\mathbb{C}}[[T]]$ (\cite[Corollary
3.3.2]{DenefLoeserBarcelona}):
\[ \mathcal{Z}_f(T) = \sum_{\emptyset\neq J\subset I} [E_J^{\circ}] \prod_{i\in J}
\frac{(\mathbb{L}-1)T^{N_i}}{\mathbb{L}^{\nu_i}-T^{N_i}}.\] In particular, this formula does
not depend on the chosen embedded resolution. For a point $x\in f^{-1}(0)$, Denef and Loeser also define the local motivic zeta function. Its formula is
\[ \mathcal{Z}_{loc,x,f}(T) = \sum_{\emptyset\neq J\subset I} [E_J^{\circ}\cap h^{-1}(x)]
\prod_{i\in J}
\frac{(\mathbb{L}-1)T^{N_i}}{\mathbb{L}^{\nu_i}-T^{N_i}}.\]

\vspace{0,4cm}

\noindent \textbf{2.2.} Let $X$ and $f$ be as above. The \textit{Hodge zeta function} of $f$ is basically the element of $\mathbb{Q}(u,v)[[T]]$ obtained by applying the Hodge-Deligne polynomial to the motivic zeta function:
\[ \mathcal{H}_f(T):= \sum_{\emptyset\neq J\subset I} H(E_J^{\circ};u,v) \prod_{i\in J} \frac{(uv-1)T^{N_i}}{(uv)^{\nu_i}-T^{N_i}}.\]
Of course one can also define the local version $\mathcal{H}_{loc,x,f}(T)$ for a point $x\in f^{-1}(0)$.\\

\noindent \textbf{2.3. Proposition.} \textsl{Let $X$ be a smooth
algebraic variety of dimension $d$ and let $f:X\to \mathbb{A}^1_{\mathbb{C}}$
be a non-constant morphism such that $X_0=f^{-1}(0)$ is irreducible, normal and
canonical (recall that a hypersurface is automatically
Gorenstein). Then \[E_{st}(X_0;u,v) = -\frac{1}{uv(uv-1)}(\mathcal{H}_f(T)(T-uv))|_{T=uv}, \] where the evaluation in $T=uv$ makes sense, since the denominator of $\mathcal{H}_f(T)$ contains the factor $T-uv$ only with multiplicity one.}\\

\noindent \textbf{Remark.} The stringy $E$-function for a hypersurface can thus be seen as a `residue' of the Hodge zeta function.\\

\noindent \textit{Proof.} Let $h:Z\to X$ be an embedded resolution
of $X_0$, with $E_i,i\in I$, the irreducible components of
$h^{-1}(X_0)$ and such that $h:h^{-1}(X_0)\to X_0$ is an isomorphism outside the singular locus $\text{Sing }X_0$ of $X_0$. For a
component $E_i$ that intersects the strict transform
$\widetilde{X_0}$ of $X_0$ (with $E_i\neq \widetilde{X_0}$) we can
look at the numerical data $(\nu_i,N_i)$ of the embedded
resolution, but also at the discrepancy $a_i$ of $\widetilde{X_0}
\cap E_i$ for $h|_{\widetilde{X_0}}: \widetilde{X_0} \to X_0$
(this is actually a log resolution). Then it is well known that
$a_i+1=\nu_i - N_i$; let us prove this for completeness. Denote by $J\subset I$ the index set of the
components $E_i\neq \widetilde{X_0}$ having nonempty intersection
with $\widetilde{X_0}$, and by $J'$ the index set of all components different from $\widetilde{X_0}$. Let $f:X_0\hookrightarrow X$ and $g:\widetilde{X_0}\hookrightarrow
Z$ be the inclusions. We have
\[ K_Z = h^{*}(K_X) + \sum_{i \in I} (\nu_i -1) E_i\]
and
\[ K_{\widetilde{X_0}}= h|_{\widetilde{X_0}}^*(K_{X_0}) + \sum_{i\in J} a_i (E_i \cap \widetilde{X_0})
\]
By the adjunction formula, the latter is also equal to
\[ g^*(K_Z + \widetilde{X_0}) = g^*(h^*(K_X) + \sum_{i\in I} (\nu_i -1)E_i) +
g^*(h^*(X_0)-\sum_{i\in J'} N_iE_i) \]
\[ = h|_{\widetilde{X_0}}^*(f^*(K_X +X_0)) + \sum_{i \in J} (\nu_i
- 1-N_i ) (E_i\cap \widetilde{X_0}),\] and then applying the adjunction
formula once more proves that $a_i+1=\nu_i - N_i$.\\

\noindent For $\widetilde{X_0}$ itself, the numerical data are
$(1,1)$. So the terms of $\mathcal{H}_f(T)$ containing a piece of the Hodge-Deligne polynomial of $\widetilde{X_0}$ assure that the
`residue' of these terms is indeed the stringy $E$-function (modulo the correction
$-\frac{1}{uv(uv-1)}$), since $X_0$ is canonical and all the
$a_i+1=\nu_i - N_i$ are thus $\geq 1$. If we can show that no
exceptional component $E_i$ with $E_i\cap X_0=\emptyset$ has
$\nu_i-N_i = 0$, then we are done. This follows from the version of inversion
of adjunction by Stevens (\cite{Stevens}). He shows that $X_0$ is canonical if and only if the pair $(X,X_0)$ is canonical near $X_0$. This means that the minimum of the $\nu_i-N_i - 1$ is greater or equal than $0$, where $i\in J'$. Alternatively, one can use the inversion of adjunction theorem for a smooth ambient variety by Ein,
Musta\c{t}\u{a} and Yasuda (see \cite[Theorem
1.6]{EinMustataYasuda}, later the first two authors generalized
this result to local complete intersection varieties in
\cite{EinMustata}). \hfill $\blacksquare$\\

\noindent \textbf{2.4. Remark.} If $X_0$ has one isolated
singular point $x$ and if we want to compute the
contribution of this singular point to the stringy $E$-function (by this we mean $E_{st}(X_0;u,v)-H(X_0\setminus\{x\};u,v)$),
we can use the formula
\[ -\frac{1}{uv(uv-1)}
(\mathcal{H}_{loc,x,f}(T)(T-uv))|_{T=uv}. \]

\section{Motivic zeta function of non-degenerate affine hypersurfaces}
\textbf{3.1.} In this section we discuss the method
of Artal, Cassou-Nogu\`es, Luengo and Melle to compute the motivic
zeta function of a polynomial that is non-degenerate with respect
to its Newton polyhedron (see \cite[Chapter 2]{Spanjaarden}). This
method is essentially earlier work by Denef and Hoornaert for Igusa's
$p$-adic zeta function (\cite{DenefHoornaert}). First we need a
lot of definitions about Newton polyhedra and polyhedral cones.
Let $f:\mathbb{A}^d_{\mathbb{C}} \to \mathbb{A}^1_{\mathbb{C}}$ be a morphism with
$f(\mathbf{0})=0$ (so $f$ is just a polynomial
$\sum_{\mathbf{n}\in (\mathbb{Z}_{\geq 0})^d}
a_{\mathbf{n}}\mathbf{x}^{\mathbf{n}}$, where
$\mathbf{x}=(x_1,\ldots,x_d),\mathbf{n}=(n_1,\ldots,n_d)$ and
$\mathbf{x}^{\mathbf{n}}=x_1^{n_1}\cdots x_d^{n_d}$). The
support of $f$ is the set $\text{supp}(f)=\{\mathbf{n}\in
(\mathbb{Z}_{\geq 0})^d \,|\, a_{\mathbf{n}}\neq 0 \}$. The Newton polyhedron
$\Gamma(f)$ of $f$ is the convex hull in $(\mathbb{R}^+)^d$ of
\[\bigcup_{\mathbf{n}\,\in\,\text{supp}(f)} \mathbf{n} +
(\mathbb{R}^+)^d.\] For the definition of a face of the Newton polyhedron we refer to \cite[p.162]{Rockafellar}. In particular, the Newton polyhedron itself is also considered as a face. A
$(d-1)$-dimensional face of the Newton polyhedron is called a
facet. For a face $\tau$ of $\Gamma(f)$, we denote
$\sum_{\mathbf{n}\in \tau} a_{\mathbf{n}} \mathbf{x}^{\mathbf{n}}$
by $f_{\tau}$. The polynomial $f$ is called \textit{non-degenerate at the origin} with
respect to its Newton polyhedron if for every compact face $\tau$
the subvariety of $(\mathbb{A}_{\mathbb{C}}^1\setminus \{0\})^d$ given by
$f_{\tau}=0$ is nonsingular. It is called \textit{non-degenerate} if the same is true for every face.\\

\noindent For $\mathbf{k}=(k_1,\ldots,k_d)\in \mathbb{R}^d$ set
$m_f(\mathbf{k}):= \inf_{\mathbf{x}\in \Gamma(f)}
\{\mathbf{k}\cdot \mathbf{x}\}$, with $\,\cdot\,$ the standard
inner product. In fact this infimum is attained and is thus a
minimum. The first meet locus of $\mathbf{k}$ is the set
$F(\mathbf{k}):= \{ \mathbf{x} \in \Gamma(f)\,|\, \mathbf{k}\cdot
\mathbf{x} = m_f(\mathbf{k}) \}$. This is a compact face of
$\Gamma(f)$ if and only if $\mathbf{k} \in (\mathbb{R}^+\setminus
\{0\})^d$. For a face $\tau$ one defines the associated polyhedral
cone $\Delta_{\tau}:=\{\mathbf{k} \in (\mathbb{R}^+)^d \,|\,
F(\mathbf{k}) = \tau\}$ in the dual space. It is well known that the cones associated with the
compact faces form a partition of $(\mathbb{R}^+\setminus
\{0\})^d$. A cone $\Delta$ is called a rational simplicial cone
(of dimension $e$) if it is generated by $e$ linearly independent
integer vectors $\mathbf{\beta}_1,\ldots,\mathbf{\beta}_e$; thus
\[ \Delta = \{ \lambda_1 \mathbf{\beta}_1 + \cdots + \lambda_e \mathbf{\beta}_e
\,|\, \lambda_i \in \mathbb{R}^+\setminus \{0\}\}.\] Usually one allows the $\lambda_i$ to be 0 in this definition, but for our goals it is more appropriate not to do that. We are interested in the set of positive integer points $\Delta'$ of such cones $\Delta$: 
\[\Delta' := \{\mathbf{\delta} \in (\mathbb{Z}_{>0})^d \,|\, n\mathbf{\delta} =
\lambda_1 \mathbf{\beta}_1 + \cdots + \lambda_e \mathbf{\beta}_e \text{ for some }
n\in \mathbb{Z}_{>0} \text{ and } \lambda_i\in \mathbb{Z}_{>0} \}.\] Let
$\mathbf{\gamma}_i$ be obtained from
$\mathbf{\beta}_i$ by dividing by the greatest common
divisor of the coordinates of $\mathbf{\beta}_i$. Then we say that
$\Delta'$ is strictly generated by
$\mathbf{\gamma}_1,\ldots , \mathbf{\gamma}_e$ and
\[ G_{\Delta'}:= \{\mathbf{\delta} \in (\mathbb{Z}_{>0})^d \,|\, \mathbf{\delta} =
\lambda_1 \mathbf{\gamma}_1 + \cdots + \lambda_e \mathbf{\gamma}_e
, 0 < \lambda_i \leq 1\}  \] is called the fundamental set of
$\Delta'$. \\

\noindent \textbf{3.2.} Every point $\mathbf{k} \in
(\mathbb{Z}_{>0})^d$ belongs to a unique cone $\Delta_{\tau}$ associated to a compact face $\tau$. Let
$\sigma(\mathbf{k})$ be $k_1 + \cdots +k_d$. Artal,
Cassou-Nogu\`es, Luengo and Melle define the following term for a
compact face $\tau$ of the Newton polyhedron of $f$ (inspired by
the work of Denef and Hoornaert):
\[S_{\Delta_{\tau}}(f, T) := \sum_{\mathbf{k} \in (\mathbb{Z}_{>0})^d
\cap \Delta_{\tau}}
\mathbb{L}^{-\sigma(\mathbf{k})}T^{m_f(\mathbf{k})}.\] Note that a priori this element does not need to belong to $\mathcal{M}_{\mathbb{C}}[[T]]$. Artal, Cassou-Nogu\`es, Luengo and Melle show that it belongs to the ring (\cite[Lemma 2.1]{Spanjaarden})
\[\mathbb{Z}[\mathbb{L},\mathbb{L}^{-1}, (1-\mathbb{L}^{-\sigma(\mathbf{a})}T^{m_f(\mathbf{a})})^{-1}][T],\] with $\mathbf{a}$ in the set of vectors such that $\mathbf{a}\,\cdot\,\mathbf{x}=M$ is a reduced integral equation of an affine hyperplane containing $\tau$. In fact they give a more general
definition for $S_{\Delta_{\tau}}$, but for our purposes this definition is sufficient.
The term can be computed by first computing a partition of
$\Delta_{\tau}$ into rational simplicial cones $\Delta_i,
i=1,\ldots, s$. Then $S_{\Delta_{\tau}}(f,T)=\sum_{i=1}^s
S_{\Delta_{i}}(f,T)$. If $\Delta_i$ is the cone strictly generated
by $\mathbf{\gamma}_1,\ldots,\mathbf{\gamma}_e$ and $G_i$ is the
fundamental set of $\Delta_i' $, then one can
prove that
\[ S_{\Delta_i}(f,T) = \left( \sum_{\mathbf{g} \in G_i }
 \mathbb{L}^{-\sigma(\mathbf{g})}T^{m_f(\mathbf{g})} \right)
 \prod_{j=1}^e
 \frac{1}{1-\mathbb{L}^{-\sigma(\mathbf{\gamma}_j)}T^{m_f(\mathbf{\gamma}_j)}}.\]

\noindent For a compact face $\tau$ Artal, Cassou-Nogu\`es, Luengo
and Melle also define a term $L_{\tau}(f,T)$, as follows (in fact they insert an extra factor $\mathbb{L}^{-d}$ in this term and in the formula for the motivic zeta function). Let
$N_{\tau}$ be the subvariety of $(\mathbb{A}^1_{\mathbb{C}}\setminus \{0\})^d $
defined by $\{f_{\tau} =0\}$ and let $[N_{\tau}]$ be its class in
the Grothendieck ring $K_0(Var_{\mathbb{C}})$. Set
\[ L_{\tau}(f,T) := (\mathbb{L}-1)^d - [N_{\tau}]+
(\mathbb{L}-1)[N_{\tau}]\frac{\mathbb{L}^{-1}T}{1-\mathbb{L}^{-1}T}
\in \mathcal{M}_{\mathbb{C}}[[T]].\] Then we can finally state the
following theorem (\cite[Theorem 2.4]{Spanjaarden}; compare with
the second remark after Theorem 4.2 in \cite{DenefHoornaert} by
`replacing' $\mathbb{L}$ by $p$ and $T$ by $p^{-s}$).\\

\noindent \textbf{Theorem.} \textsl{Let $f$ be a polynomial in $d$
variables over the complex numbers with $f(\mathbf{0})=0$. Assume
that $f$ is non-degenerate at the origin with respect to its Newton polyhedron
$\Gamma(f)$. Then
\[ \mathcal{Z}_{loc,0,f}(T)= \sum_{\substack{\text{compact faces}\\ \tau\text{ of }\Gamma(f)}} L_{\tau}(f,T)
S_{\Delta_{\tau}}(f,T).\]}

\noindent When $f$ is non-degenerate with respect to $\Gamma(f)$, an analogous formula for $\mathcal{Z}_f(T)$ could be given, by
summing over \textit{all} faces of $\Gamma(f)$ (see Theorem 4.2 in
\cite{DenefHoornaert}). As an easy corollary of this theorem, of Remark 2.4 and of \cite[Theorem (4.6)]{YPG} we obtain:\\

\noindent \textbf{Corollary.} \textsl{Let $f$ be a polynomial in $d$
variables over the complex numbers with $f(\mathbf{0})=0$. Assume
that $f$ is non-degenerate at the origin with respect to its Newton polyhedron
$\Gamma(f)$ and that $f^{-1}(0)$ has only an isolated canonical singularity at the origin. Then the contribution of this singular point to the stringy $E$-function is given by
\[ \sum_{\substack{\text{compact faces}\\ \tau\text{ of }\Gamma(f)}} H(N_{\tau};u,v) \tilde{S}_{\Delta_{\tau}}(f,uv)  ,\]
where $\tilde{S}_{\Delta_{\tau}}(f,uv)$ is obtained from $S_{\Delta_{\tau}}(f,T)$ by `replacing' both $\mathbb{L}$ and $T$ by $uv$.}

\section{Stringy $E$-functions of Brieskorn singularities}
\noindent \textbf{4.1.} Using the results described in the previous sections, we want to sketch an easy way to compute the contribution of a Brieskorn singularity to the stringy $E$-function. A Brieskorn singularity is given by the origin of the zero set in $\mathbb{A}_{\mathbb{C}}^d$ of a polynomial of the form
\[f(x_1,\ldots,x_d) = x_1^{a_1}+\cdots +x_d^{a_d},\] where all $a_i\geq 2$. Put $k:=\text{lcm}(a_1,\ldots,a_d)$, $\alpha:=(\frac{k}{a_1},\ldots,\frac{k}{a_d})\in \mathbb{Z}_{>0}^d$ and $\Sigma:=\sum_i \frac{k}{a_i}$. 
To describe when Brieskorn singularities are canonical or terminal, we use Proposition (4.3) from \cite{Reid2} and Theorem (4.6) of \cite{YPG} applied to $\alpha$. We get the following:
\begin{itemize}
\item[(1)] the above singularity is canonical if and only if $\Sigma - k \geq 1$, 
\item[(2)] if the above singularity is terminal then $\Sigma - k \geq 2$.
\end{itemize}
In particular, if $\Sigma - k =1$, then the singularity is strictly canonical. From now on, we assume that $\Sigma - k \geq 1$.\\

\noindent \textbf{Remark.} The criterion for terminality given in \cite{Lin} is not correct, as can be seen for instance from the $A_n$ surface singularities for even $n$.\\

\noindent Denote by $I$ the index set $\{1,\ldots, d\}$. Let $\mathcal{S}$ be the set of subsets $J\varsubsetneq I$ defined by
\[ J\in \mathcal{S} \Leftrightarrow \begin{array}{l} J=\emptyset \text{ or for all }j'\in J \text{ we have that }\\ \gcd\{\alpha_j\,|\, j\in I\setminus J\} > \gcd\{\alpha_{j}\,|\,j \in \{j'\}\cup (I\setminus J)\}.\end{array}\]
For example, if $\alpha = (6,6,4,3,3)$ then $\mathcal{S}=\{\emptyset,\{3\},\{4,5\},\{3,4,5\},\{1,2,4,5\}\}$. For any $J\varsubsetneq I$, denote $\gcd\{\alpha_j\,|\, j\in I\setminus J\}$ by $g_J$. It is easy to see that all cones associated to compact faces of the Newton polyhedron of $f$ are strictly generated by $\alpha$ and between 0 and $d-1$ standard basis vectors $\mathbf{e}_i$. To $J\varsubsetneq I$ we associate the cone $\Delta_J$ generated by $\alpha$ and the $\mathbf{e}_j$ for $j \in J$. The compact face corresponding to this cone is denoted by $\tau_J$. The fundamental set $G_J$ of such a cone is given by
\[ \left\{\delta_J^l:=\frac{l}{g_{J}} \alpha + \sum_{j\in J} \frac{g_{J} -  (l \alpha_j\ \text{mod}\, g_{J})}{g_{J}} \mathbf{e}_j\,|\, l=1,\ldots , g_{J} \right\}.\]
Note that we need the numbers $\sigma(\delta_J^l)$ and $m_f(\delta_J^l)=\frac{kl}{g_J}$ for the formula of the stringy $E$-function (see Section 3). A short computation shows that
\[ \Sigma - k + |J| - \sigma(\delta_J^l) + m_f(\delta_J^l) \geq 0 \] for all $J$ and $l$. We also have the following lemma.\\

\noindent \textbf{Lemma.} \textsl{Let $J \varsubsetneq I$.
\begin{itemize}
\item[(1)] For $J'\subset J$
\[ (uv)^{\Sigma - k + |J|}\biggl( \sum_{l=1}^{g_J} (uv)^{-\sigma(\delta_J^l)+m_f(\delta_J^l)} \biggr) - (uv)^{\Sigma - k + |J'|}\biggl( \sum_{l'=1}^{g_{J'}} (uv)^{-\sigma(\delta_{J'}^{l'})+m_f(\delta_{J'}^{l'})} \biggr) \]
is a polynomial in $uv$ with nonnegative coefficients.
\item[(2)] $J\in \mathcal{S}$ if and only if $J=\emptyset$ or the above polynomial is nonzero for all $J' \varsubsetneq J$.
\end{itemize}}

\noindent \textit{Proof.} If $J' \subset J$ then $g_{J'}\,|\,g_J$. For $l'\in \{1,\ldots,g_{J'}\}$ we take $l=\frac{l'g_J}{g_{J'}}$. Then 
\[ \delta_J^l = \delta_{J'}^{l'} + \sum_{j\in J\setminus J'} \mathbf{e}_j\]
and (1) follows. For (2) we note that $\emptyset \neq J \in \mathcal{S}$ if and only if $g_{J'}\neq g_J$ for all $J' \varsubsetneq J$. And then the proof of (1) implies (2). \hfill $\blacksquare$\\

\noindent \textbf{4.2.} For $J\varsubsetneq I$ we can now define a polynomial $p_J(uv)$ in $uv$ in a recursive way as follows. For $J= \emptyset$ it is just 1. Note that this equals \[(uv)^{\Sigma - k + |J|}\biggl( \sum_{l=1}^{g_J} (uv)^{-\sigma(\delta_J^l)+m_f(\delta_J^l)} \biggr).\] For $J\neq \emptyset$ we define $p_J$ as 
\[(uv)^{\Sigma - k + |J|}\biggl( \sum_{l=1}^{g_J} (uv)^{-\sigma(\delta_J^l)+m_f(\delta_J^l)} \biggr) - \sum_{J' \varsubsetneq J} p_{J'}.\]

\vspace{0.3cm}

\noindent \textbf{Lemma.} \textsl{The polynomial $p_J$ is nonzero if and only if $J\in \mathcal{S}$, and in that case it has nonnegative coefficients.}\\

\noindent \textit{Proof.} From the proof of Lemma 4.1 it is clear that the vectors $\delta_J^l$ that cannot be written as
\[ \delta_{J'}^{l'} + \sum_{j\in J\setminus J'} \mathbf{e}_j \] for $J'\varsubsetneq J$ (and necessarily $l'=\frac{lg_{J'}}{g_J}$) give a contribution to $p_J$ with nonnegative coefficients. Moreover, there are such vectors if and only if $J\in \mathcal{S}$. Now we only have to show that we subtract the contribution of the vectors $\delta_J^l$ that can be written as $\delta_{J'}^{l'} + \sum_{j\in J\setminus J'} \mathbf{e}_j$ exactly once. For such a $\delta_J^l$, the complement $J'$ of the set of $j$'s for which the coefficient $\frac{g_{J} -  (l\alpha_j\ \text{mod}\,g_{J})}{g_{J}}$ of $\mathbf{e}_j$ in the definition of $\delta_J^l$ is equal to 1, is the minimal set for which $\delta_J^l$ can be written as $\delta_{J'}^{l'} + \sum_{j\in J\setminus J'} \mathbf{e}_j$, for $l'=\frac{lg_{J'}}{g_J}$. So the contribution 
\[(uv)^{\Sigma - k + |J|}(uv)^{-\sigma(\delta_J^l)+m_f(\delta_J^l)} =  (uv)^{\Sigma - k + |J'|}(uv)^{-\sigma(\delta_{J'}^{l'})+m_f(\delta_{J'}^{l'})}\] is subtracted exactly in $p_{J'}$.\hfill $\blacksquare$\\

\noindent For a compact face $\tau$, the polynomial $f_{\tau}$ contains in this case only $\dim\tau +1$ variables; so we can consider its zero set in $\mathbb{A}_{\mathbb{C}}^{\dim\tau +1}$. We denote this zero set by $M_{\tau}$. Finally we have settled all notations for the following theorem.\\

\noindent \textbf{Theorem.} \textsl{The contribution of a Brieskorn singularity to the stringy $E$-function is given by
\[ \frac{1}{(uv)^{\Sigma-k}-1} \left(\sum_{J\in \mathcal{S}} \bigl( H(M_{\tau_J};u,v) -1\bigr) p_J(uv)\right).\]}

\noindent \textit{Proof.} For a compact face $\tau$ of the Newton polyhedron of $f$ we have considered the zero set of $f_{\tau}$ in Section 3 in $(\mathbb{A}_{\mathbb{C}}^1 \setminus \{0\})^d$, and we have denoted this zero set by $N_{\tau}$. Since $f_{\tau}$ has only $\dim \tau +1$ variables, we can as well consider the zero set of $f_{\tau}$ in $(\mathbb{A}_{\mathbb{C}}^1 \setminus \{0\})^{\dim \tau +1}$, and we call this zero set $\widetilde{N_{\tau}}$. In particular, 
\[H(N_{\tau};u,v) = (uv-1)^{d-\dim\tau -1} H(\widetilde{N_{\tau}};u,v).\] According to Corollary 3.2 we have to compute
\[ A := \sum_{\substack{\text{compact faces}\\ \tau\text{ of }\Gamma(f)}} H(N_{\tau};u,v) \tilde{S}_{\Delta_{\tau}}(f,uv)  ,\]
which in this case can be written as
\[\sum_{J\varsubsetneq I} H(N_{\tau_J};u,v) \tilde{S}_{\Delta_J}(f,uv).\]
This equals
\[ \frac{(uv)^{\Sigma - k}}{(uv)^{\Sigma - k}-1} \left(\sum_{J\varsubsetneq I} (uv)^{|J|}  H(\widetilde{N_{\tau_J}};u,v) \biggl(\sum_{\mathbf{g}\in G_J} (uv)^{-\sigma(\mathbf{g})+m_f(\mathbf{g})} \biggr)\right) . \tag{1} \] 
We also have
\begin{eqnarray*}
H(\widetilde{N_{\tau_J}};u,v) & = &\left( \sum_{I\varsupsetneq J' \supseteq J} (-1)^{|J'|-|J|} H(M_{\tau_{J'}};u,v) \right) + (-1)^{d-|J|} \\
 & = & \sum_{I\varsupsetneq J' \supseteq J} (-1)^{|J'|-|J|} (H(M_{\tau_{J'}};u,v) -1).
\end{eqnarray*}
Putting this into (1) and summing over $J'\varsubsetneq I$ leads to 
\begin{eqnarray*}
A &=& \frac{1}{(uv)^{\Sigma-k}-1} \left(\sum_{J'\varsubsetneq I} \bigl( H(M_{\tau_{J'}};u,v) -1\bigr) p_{J'}(uv)\right)\\
 &=&  \frac{1}{(uv)^{\Sigma-k}-1} \left(\sum_{J'\in \mathcal{S}} \bigl( H(M_{\tau_{J'}};u,v) -1\bigr) p_{J'}(uv)\right). 
\end{eqnarray*}
\hfill $\blacksquare$\\

\noindent \textbf{4.3.} To apply the above theorem in concrete examples, we only have to explain how one can compute the Hodge-Deligne polynomial of an $M_{\tau_{J}}$. The equation of an $M_{\tau_{J}}$ is trivially quasi-homogeneous and thus it can be computed by the following method described in \cite[Section 2]{Dais}. Let $f\in \mathbb{C}[x_1,\ldots,x_{r+1}]$ be a quasihomogeneous polynomial of degree $d$ with respect to the weights $w_1,\ldots,w_{r+1}$ and assume that $\mathbf{0}$ is an isolated singularity of $Y:=f^{-1}(0)$. It is well known that the singular cohomology of the link $L$ of the singularity carries a natural mixed Hodge structure. According to \cite[Proposition 2.8]{Dais},
\[H(Y;u,v) = (uv)^r +(-1)^{r-1}(uv-1)\sum_{p=0}^{r-1} h^{p,r-1-p}(H^{r-1}(L,\mathbb{C}))u^pv^{r-1-p}, \]
where $h^{p,r-1-p}(H^{r-1}(L,\mathbb{C}))$ denotes the dimension of the $H^{p,r-1-p}$ component of the mixed Hodge structure on the cohomology group $H^{r-1}(L,\mathbb{C})$. These numbers can be computed in terms of $w_1,\ldots,w_{r+1}$, as explained in Theorem 2.6 and Lemma 2.7 of \cite{Dais}. Consider the Milnor algebra 
\[ M(f):=\frac{\mathbb{C}[x_1,\ldots,x_{r+1}]}{\left(\frac{\partial f}{\partial x_1},\ldots,\frac{\partial f}{\partial x_{r+1}} \right)}.      \]
This becomes a finitely generated graded $\mathbb{C}$-algebra if we give $x_i$ degree $w_i$. The Poincar\'e series of such an algebra is defined by \[ P_{M(f)}(t) := \sum_{k\geq 0} (\dim_{\mathbb{C}}M(f)_k)t^k,     \] where $M(f)_k$ is the piece of degree $k$. This series can be calculated by the formula
\[ P_{M(f)}(t)= \frac{(1-t^{d-w_1})\cdots (1-t^{d-w_{r+1}})}{(1-t^{w_1})\cdots (1-t^{w_{r+1}})}.  \]
Dais shows, referring to work of Griffiths and Steenbrink (\cite{Griffiths} and \cite{Steenbrinkquasihomogeen}), that the numbers $h^{p,r-1-p}(H^{r-1}(L,\mathbb{C}))$ equal \[\dim_{\mathbb{C}}M(f)_{(p+1)d-(w_1+\cdots+w_{r+1})},\] and thus they can be computed from the Poincar\'e series.\\

\noindent \textbf{4.4.} \textbf{Corollary.} \textsl{The contribution of a Brieskorn singularity to the stringy $E$-function can be written in the form
\[ \frac{P(u,v)} {(uv)^{\Sigma - k - 1} + (uv)^{\Sigma - k -2} + \cdots + 1},\]
where $P(u,v)$ is a polynomial $\sum_{i,j} c_{i,j} u^iv^j$ whose coefficients satisfy $(-1)^{i+j} c_{i,j} \geq 0$. In particular, if $\Sigma - k =1$ then this contribution is a polynomial. Furthermore, if $Y$ is a projective variety with at most Brieskorn singularities with $\Sigma - k=1$, then its stringy Hodge numbers are nonnegative.}\\

\noindent \textit{Proof.} For the first statement we combine Lemma 4.2, Theorem 4.2 and the discussion in 4.3. Indeed, it is clear that $H(M_{\tau_{J}};u,v)-1$ is divisible by $uv-1$ and that this division leads to a numerator of the requested form. The last statement of the corollary follows then trivially from the more general lemma below. \hfill $\blacksquare$ \\

\noindent \textbf{Lemma.} \textsl{Let $Y$ be a projective variety with at most isolated Gorenstein canonical singularities and with a polynomial stringy $E$-function. Assume that the contribution of the singularities to the stringy $E$-function is a polynomial $\sum_{i,j} c_{i,j} u^iv^j$ whose signs are `right' in the sense that $(-1)^{i+j}c_{i,j} \geq 0$. Then the stringy Hodge numbers of $Y$ are nonnegative.}\\

\noindent \textit{Proof.} By Theorem (1.13) from \cite{Steenbrinkpuriteit} $H^i(Y)=H^i_c(Y)$ has a pure Hodge structure of weight $i$ for $i > d$, where $d$ is the dimension of $Y$. So $h^{p,q}(H^i(Y)) = 0$ for $i>d, p+q \neq i$. Since $Y$ is projective, $h^{p,q}(H^i(Y))=0$ for $p+q > i$, where $i$ is now arbitrary (see \cite[Th\'eor\`eme (8.2.4)]{Deligne2}). So if we fix $(p,q)$ with $p+q \geq d$, the coefficient of $u^pv^q$ in $H(Y;u,v)$ is exactly 
\[ (-1)^{p+q} h^{p,q}(H^{p+q}(Y)) \]
and thus it has the `right' sign. This is also the coefficient of $u^pv^q$ in $H(Y_{ns};u,v)$, where $Y_{ns}$ denotes the nonsingular part of $Y$, because $H(Y;u,v)-H(Y_{ns};u,v)$ equals the number of singular points of $Y$. So the coefficients of $u^pv^q$ in $E_{st}(Y)$ also have the `right' sign for $p+q \geq d$ and by symmetry this is then true for all coefficients. \hfill $\blacksquare$\\

\noindent \textbf{4.5.} \textbf{Remark.} 
\begin{itemize}
\item[(1)] We think that it is very difficult to give a combinatorial description of the Brieskorn singularities with $\Sigma - k \geq 2$ that give a polynomial contribution to the stringy $E$-function. To get an idea of what one can expect we refer to \cite[Theorem A]{Lin}, where for some specific Brieskorn singularities the existence of a crepant resolution is studied. In these cases, the singularities that give a polynomial contribution to the stringy $E$-function seem to correspond to the ones admitting a crepant resolution (in general there exist Brieskorn singularities with a polynomial contribution that do not admit a crepant resolution, see the example in \cite{Schepers}). 
\item[(2)] Let us compare Corollary 4.4 with the theorem of \cite{SchepersVeys} discussed at the end of 1.6 and in 5.1 below. From the point of view of that theorem, the strictly canonical singularities are the `worst' ones. So it is somewhat surprising that we obtain here the nonnegativity of the stringy Hodge numbers exactly for a class of strictly canonical Brieskorn singularities. We do not know how to prove it for Brieskorn singularities with $\Sigma-k \geq 2$ and a polynomial contribution to the stringy $E$-function. 
\end{itemize}

\section{An interesting example}
\noindent \textbf{5.1.} We want to conclude this paper by computing a concrete stringy $E$-function of a projective variety with Brieskorn singularities. Before we give the details of the example, let us first explain why it is interesting. As already mentioned in the introduction, we proved the following theorem in \cite{SchepersVeys}.\\

\noindent \textbf{Theorem.} \textsl{Let $Y$ be a $d$-dimen\-sional
Gorenstein projective variety with at most isolated singularities ($d\geq 3$). Let $f:X\to Y$ be a log resolution. 
Assume that the discrepancy coefficients of the exceptional components are strictly greater than $\lfloor
\frac{d-4}{2} \rfloor$ (this condition does not depend on the chosen log resolution). Write the stringy $E$-function of $Y$ as a power series $\sum_{i,j\geq 0}b_{i,j}u^iv^j$. Then $(-1)^{i+j}b_{i,j}\geq 0$ for $i+j\leq d$. In particular, if the stringy $E$-function of $Y$ is a polynomial, then Ba\-ty\-rev's conjecture is true for $Y$. For $d=3$ the statements remain true if we drop the hypothesis of isolated singularities.}\\

\noindent Our example shows that this theorem cannot be extended to the case of 6-dimensional varieties with terminal singularities (so now we also allow discrepancy coefficients equal to 1). It is an example of a non-polynomial stringy $E$-function with a negative number $b_{3,3}$.\\

\noindent \textbf{5.2.} For the example we also need the formula for the Hodge-Deligne polynomial of a Fermat hypersurface. This is explained in \cite{Dais}. We denote the $d$-dimensional Fermat hypersurface of degree $l$ by $Y^{(d)}_l$. So $Y^{(d)}_l$ is given by 
\[ \{x_0^l + \cdots + x_{d+1}^l =0\} \subset \mathbb{P}^{d+1}_{\mathbb{C}}.    \]
To write down the Hodge-Deligne polynomial of $Y^{(d)}_l$ we need an auxiliary definition. Dais considers the numbers 
\[ \mathcal{G}(\kappa,\lambda\,|\,\nu,\xi) := \sum_{j=0}^{\lambda} (-1)^j {\kappa + 1 \choose j} {\nu(\lambda -j) + \xi \choose \kappa} \] for $(\kappa,\lambda,\nu,\xi)\in \mathbb{Z}_{\geq 0}^4$ and $\kappa \geq \lambda$ (if $m>n$, the binomial coefficient ${n\choose m}$ must be interpreted as 0). Then the Hodge-Deligne polynomial of $Y^{(d)}_l$ is given by (\cite[Lemma 3.3]{Dais})
\[ H(Y^{(d)}_l;u,v) := \sum_{p=0}^{d} u^p\left(v^p + (-1)^d\mathcal{G}(d+1,p+1\,|\,l-1,p)v^{d-p} \right). \]

\vspace{0.3cm}

\noindent \textbf{5.3.} \textit{Computation of the example using Theorem 4.2.} We want to compute the stringy $E$-function of the singular variety
\[ Y := \{ x_1^5z + x_2^5z +x_3^6+x_4^6+x_5^6+x_6^6+x_7^6 =0 \} \subset \mathbb{P}^7_{\mathbb{C}}, \]
where we consider $z=0$ as the hyperplane at infinity. There are 6 isolated singularities: the origin of the affine chart $z\neq 0$ (local equation $x_1^5+x_2^5+x_3^6+x_4^6+x_5^6+x_6^6+x_7^6 =0$) and five singularities at infinity, all analytically isomorphic to the origin of $\{x_1^2+x_2^2+x_3^6+x_4^6+x_5^6+x_6^6+x_7^6 =0\}$. Let us start with the first singularity. Using the notations from Section 4 we have
\begin{itemize}
\item $k=30, \alpha = (6,6,5,5,5,5,5), \Sigma = 37$,
\item $I=\{1,\ldots,7\},\mathcal{S}=\{\emptyset,\{1,2\},\{3,4,5,6,7\}\}$.
\end{itemize}
An easy computation shows that
\[\begin{array}{l}
p_{\emptyset} = 1,\\
p_{\{1,2\}} = (uv)^6 + (uv)^5 + (uv)^4 + (uv)^3, \\ 
p_{\{3,4,5,6,7\}} = (uv)^{10} + (uv)^8 + (uv)^6 + (uv)^4 + (uv)^2.
\end{array}\]
Using the discussion in 4.3 one computes
\[\begin{array}{l}
H(M_{\tau_{\emptyset}};u,v) = (uv)^6-(uv-1)(20u^4v+20uv^4+1020u^3v^2+1020u^2v^3),\\
H(M_{\tau_{\{1,2\}}};u,v) = (uv)^4 -(uv-1)(5u^3 +5v^3 +255u^2v+255uv^2), \\ 
H(M_{\tau_{\{3,4,5,6,7\}}};u,v) = 5uv-4.
\end{array}\]
In this way the contribution of this singularity to the stringy $E$-function becomes
\begin{equation*} \begin{split}
A:= & \frac{(uv-1)}{((uv)^7-1)}\biggl( 5(uv)^{10}+(uv)^9+7(uv)^8+3(uv)^7+9(uv)^6+4(uv)^5 +8(uv)^4\\ & +2(uv)^3  +6(uv)^2+uv+1 -5u^9v^6-5u^6v^9 -255u^8v^7   -255v^7v^8-5u^8v^5\\ & -5u^5v^8 -255u^7v^6 -255u^6v^7-5u^7v^4-5u^4v^7    -255u^6v^5  -
255u^5v^6 -5u^6v^3\\ & -5u^3v^6-255u^5v^4 -255u^4v^5 -20u^4v-20uv^4-1020u^3v^2-1020u^2v^3\biggr).
\end{split} \end{equation*}
The contribution of a singularity at infinity is easier to compute and equals
\begin{equation*} \begin{split} B:=\frac{(uv-1)}{((uv)^5-1)} & \bigg(5(uv)^5+(uv)^4+(uv)^3+(uv)^2+uv+1 \\ & -5u^4v-5uv^4 -255u^3v^2-
255u^2v^3\bigg). \end{split}  \end{equation*}
Now we still have to compute the contribution of the nonsingular part of $Y$. Let us first do that at infinity. The total part at infinity is given by
\[ Y^{\infty}:=\{x_3^6+x_4^6+x_5^6+x_6^6+x_7^6=0 \} \subset \mathbb{P}^6_{\mathbb{C}}.\]
To find the nonsingular part, we just have to remove five points. In fact $Y^{\infty}$ can be found from the Fermat hypersurface $Y^{(3)}_6$ (notation as in 5.2) by taking twice the projective cone. On the level of the Hodge-Deligne polynomial, one such operation multiplies the original Hodge-Deligne polynomial by $uv$ and adds 1. By 5.2 \[H(Y_6^{(3)};u,v)= (uv)^3+(uv)^2+uv+1-5u^3-5v^3-255u^2v-255uv^2.\] The contribution at infinity becomes thus (do not forget to subtract the 5 singular points)
\[ C := (uv)^5+(uv)^4+(uv)^3+(uv)^2+uv -4 -5u^5v^2-5u^2v^5-255u^4v^3-255u^3v^4.\]
Finally, we only have to compute the contribution of $\{x_1^5+x_2^5+x_3^6+\cdots +x_7^6=0\}\subset \mathbb{A}_{\mathbb{C}}^7$ minus the singular point. This can be done by the method of 4.3 and the result is (subtracting also the singular point)
\[ D := (uv)^6-1-(uv-1)(20u^4v+20uv^4+1020u^3v^2+1020u^2v^3).\]
To find the total stringy $E$-function of $Y$ we just add $A,5B,C$ and $D$ and simplify. The result is
\begin{equation*} \begin{split} & E_{st}(Y;u,v) =  \frac{1}{((uv)^5-1)((uv)^7-1)}  \biggl( (uv)^{18}+(uv)^{17}+6(uv)^{16}-3(uv)^{15}\\ & +7(uv)^{14} +21(uv)^{13}-20(uv)^{12}-12(uv)^{11}+6(uv)^{10}-14(uv)^9 
+6(uv)^8\\ &  -12(uv)^7 -20(uv)^6+21(uv)^5+7(uv)^4-3(uv)^3+6(uv)^2+uv+1\\ &  -25(u^{17}v^{14} +u^{14}v^{17}+u^4v+uv^4) -1275(u^{16}v^{15}+u^{15}v^{16}+u^3v^2+u^2v^3)\\ & +20(u^{16}v^{13} +u^{13}v^{16}+u^5v^2+u^2v^5)+
1020(u^{15}v^{14}+u^{14}v^{15}+u^4v^3+u^3v^4)\\ & -5(u^{15}v^{12} +u^{12}v^{15}+u^6v^3+u^3v^6) -255(u^{14}v^{13}+u^{13}v^{14}+u^5v^4+u^4v^5)\\ &  +10(u^{11}v^8 +u^8v^{11}+u^{10}v^7+u^7v^{10}) + 510(u^{10}v^9+u^9v^{10}+u^9v^8+u^8v^9)      \biggr).  \end{split} \end{equation*}
So if we develop this in power series, we get a term $-3(uv)^3$. To understand why the Theorem 5.1 does not extend to this case, it is useful to compute an explicit log resolution, although this is a bit complicated. \\

\noindent \textbf{5.4.} \textit{Computation of the example using a log resolution.} Let us consider the singularities at infinity first. These are all analytically isomorphic to the origin of 
\[ Y':=\{x_1^2+x_2^2+x_3^6+\cdots+x_7^6 = 0 \} \subset \mathbb{A}^7_{\mathbb{C}}. \] To resolve them, we have to blow up in the singular point first. This gives two exceptional components, denoted $D_1^{\infty}$ and $D_2^{\infty}$, and their intersection becomes the new singular locus. Blowing it up, gives again two exceptional components $E_1^{\infty}$ and $E_2^{\infty}$ whose intersection is singular for the strict transform of $Y'$. Moreover, the components $D_1^{\infty}$ and $D_2^{\infty}$ are separated (we use the same name for a divisor at any stage of the resolution process, instead of speaking of the strict transform). Then we only have to perform the blow-up in the intersection of $E_1^{\infty}$ and $E_2^{\infty}$. This gives one new exceptional component $F^{\infty}$ and the following intersection diagram:
\begin{center}
\begin{picture}(65,10)
\put(2,2){\circle*{2}} \put(17,2){\circle*{2}}
\put(32,2){\circle*{2}} \put(47,2){\circle*{2}}
\put(62,2){\circle*{2}} \put(2,2){\line(1,0){60}}
\put(0,6){$D_1^{\infty}$} \put(15,6){$E_1^{\infty}$}
\put(30,6){$F^{\infty}$} \put(45,6){$E_2^{\infty}$} \put(60,6){$D_2^{\infty}$}
\end{picture}
\end{center}
The discrepancy coefficient of all these components is 4. The components $D_1^{\infty}$ and $D_2^{\infty}$ are isomorphic to $\mathbb{P}^5_{\mathbb{C}}$, $E_1^{\infty}$ and $E_2^{\infty}$ are $\mathbb{P}^1_{\mathbb{C}}$-bundles over $\mathbb{P}^4_{\mathbb{C}}$ and all intersections are isomorphic to $\mathbb{P}^4_{\mathbb{C}}$. It is not so easy to compute the Hodge-Deligne polynomial of $F^{\infty}$. In one of the charts, $F^{\infty}$ is given by the equations 
\[ \{ x_3=0 , x_1^2 +x_2^2 + 1 +x_4^6+x_5^6+x_6^6+x_7^6 = 0\} \subset \mathbb{A}^7_{\mathbb{C}}.\]
This variety is isomorphic to 
\[ \{ x_1x_2 + 1 +x_4^6 +\cdots +x_7^6 =0\} \subset \mathbb{A}^6_{\mathbb{C}}.\]
For $x_1\neq 0$, one finds a contribution of $(uv-1)(uv)^4$ to $H(F^{\infty};u,v)$. For $x_1 = 0$, one finds $uv$ times the Hodge-Deligne polynomial of an affine piece of the 3-dimensional Fermat hypersurface $Y_6^{(3)}$ of degree 6 (in fact $Y_6^{(3)}\setminus Y_6^{(2)}$), and this Hodge-Deligne polynomial can be calculated by 5.2. Taking into account the contributions of all other relevant coordinate charts (that can be calculated analogously), one finds 
\[ H(F^{\infty};u,v) = (uv)^5+2(uv)^4+2(uv)^3+2(uv)^2+2(uv)+1+(uv)H(Y^{(3)}_6;u,v),\] with $H(Y^{(3)}_6;u,v)$ as above. With these data, one can compute that the contribution of such a singular point to the stringy $E$-function is indeed formula $B$ from 5.3.\\

\noindent The computation of the contribution of the singularity in the origin of the affine chart $z\neq 0$ of $Y$ can be done as follows. First we blow up in the singular point itself. This produces five exceptional components $D_1,\ldots,D_5$. After this blow up, they all go through one $\mathbb{P}^4_{\mathbb{C}}$, and thus they have nowhere normal crossings. The new singular locus is a Fermat hypersurface of degree 6 on this $\mathbb{P}^4_{\mathbb{C}}$. Blowing up in this singular locus gives two new exceptional divisors $E_1$ and $E_2$. They both intersect the $D_i$, but the new singular locus is the intersection of $E_1$ and $E_2$. This singular locus also contains a piece of the intersection of the $D_i$, and this piece is exactly the intersection of the $D_i$ with $E_2$ (which is surprisingly only 3-dimensional). Blowing up in this new singular locus splits of $E_2$ from $E_1$ and all the $D_i$. Two new exceptional components $F_1$ and $F_2$ appear. They intersect each other, and apart from that, the first intersects $E_1$ and all the $D_i$, and the second intersects $E_2$ and the $D_i$. The new singular locus consists of five separate pieces; one piece on each $D_i$. It is exactly the intersection of $F_2$ with the $D_i$. Blowing it up gives us five new components $G_1,\ldots, G_5$, all intersecting $F_1,F_2$, and every $G_i$ intersects one $D_i$ (of course we take a compatible numbering). Finally we have a nonsingular strict transform, but still the $D_i$ have no normal crossings. Blowing up in their intersection (which is isomorphic to $\mathbb{P}^4_{\mathbb{C}}$), gives one new exceptional component $C$ (a $\mathbb{P}^1_{\mathbb{C}}$-bundle over this intersection), intersecting each $D_i$ and also intersecting $E_1$. We find the following intersection diagram:
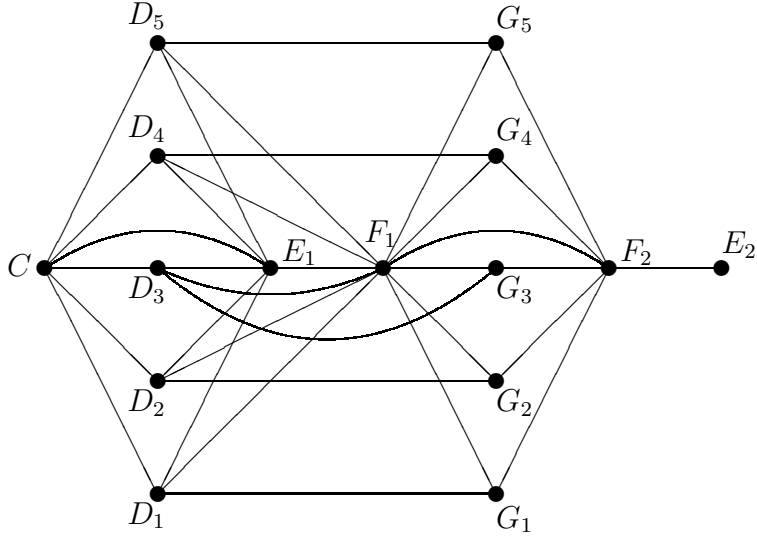
\begin{center}
\begin{picture}(102,70)
\put(20,5){\circle*{2}} \put(20,20){\circle*{2}} \put(20,35){\circle*{2}} \put(20,50){\circle*{2}} \put(20,65){\circle*{2}} \put(5,35){\circle*{2}} \put(35,35){\circle*{2}} \put(50,35){\circle*{2}} \put(65,35){\circle*{2}} \put(65,5){\circle*{2}} \put(65,20){\circle*{2}} \put(65,50){\circle*{2}} \put(65,65){\circle*{2}} \put(80,35){\circle*{2}} \put(95,35){\circle*{2}} \put(5,35){\line(1,0){90}} \put(5,35){\line(1,1){15}} \put(5,35){\line(1,2){15}} \put(5,35){\line(1,-1){15}} \put(5,35){\line(1,-2){15}} \put(35,35){\line(-1,-1){15}} \put(35,35){\line(-1,-2){15}} \put(35,35){\line(-1,1){15}} \put(35,35){\line(-1,2){15}} \put(50,35){\line(-2,1){30}} \put(20,65){\line(1,-1){30}} \put(50,35){\line(-2,-1){30}} \put(20,5){\line(1,1){30}} \put(50,35){\line(1,1){15}} \put(50,35){\line(1,2){15}} \put(50,35){\line(1,-1){15}} \put(50,35){\line(1,-2){15}} \put(80,35){\line(-1,1){15}} \put(80,35){\line(-1,2){15}} \put(80,35){\line(-1,-1){15}} \put(80,35){\line(-1,-2){15}} \qbezier(5,35)(20,45)(35,35) \qbezier(50,35)(65,45)(80,35) \put(20,65){\line(1,0){45}} \put(20,50){\line(1,0){45}} \put(20,20){\line(1,0){45}} \put(20,5){\line(1,0){45}} \qbezier(20,35)(42.5,16)(65,35) \qbezier(20,35)(35,28)(50,35) \put(0,34){$C$} \put(16,1){$D_1$} \put(16,16){$D_2$} \put(16,31){$D_3$} \put(16,52.5){$D_4$} \put(16,67.5){$D_5$} \put(36.5,36){$E_1$} \put(47.5,38.5){$F_1$} \put(65,0.5){$G_1$} \put(65,16){$G_2$} \put(65,31){$G_3$} \put(65,52){$G_4$} \put(65,67){$G_5$} \put(81.5,36){$F_2$} \put(95,37){$E_2$}
\end{picture}
\end{center}
The discrepancy coefficients of the $D_i$ are 1. One finds 6 for $C$, 5 for $E_1$ and the $G_i$, 4 for $F_1$, 3 for $F_2$ and 2 for $E_2$. There are twenty threefold intersections, namely $C\cap E_1 \cap D_i$, $E_1\cap F_1\cap D_i$, $F_1 \cap D_i \cap G_i$ and $F_1\cap F_2 \cap G_i$, where $i$ runs of course from 1 to 5. They are all isomorphic to the Fermat hypersurface $Y^{(3)}_6$. One can count from the diagram that there are thirty-four twofold intersections, all having Hodge-Deligne polynomial $(uv+1)H(Y_6^{(3)};u,v)$, except for the $C\cap D_i$, they are isomorphic to $\mathbb{P}^4$. The Hodge-Deligne polynomials of the components itself are
\[ \begin{array}{l}
H(C)=(uv+1)((uv)^4+(uv)^3+(uv)^2+uv+1),\\
H(D_i)=(uv)^5+(uv)^4+(uv)^3+(uv)^2+uv+1 + 3uvH(Y_6^{(3)};u,v), \\
H(E_1)=H(G_i)=((uv)^2+2uv+1)H(Y_6^{(3)};u,v),\\
H(E_2)=((uv)^2+uv+1)H(Y_6^{(3)};u,v),\\
H(F_1)=H(F_2)=((uv)^2+7uv+1)H(Y_6^{(3)};u,v).
\end{array} \]
By a rather lengthy calculation one can then simplify the contribution of this singular point indeed to expression $A$ from 5.3.\\

\noindent To conclude this example, we have a closer look at the coefficient $b_{3,3}$ from the power series development $\sum_{i,j\geq 0} b_{i,j}u^iv^j$ of the stringy $E$-function. Let $Z$ be a general Gorenstein canonical projective variety of dimension 6. Take a log resolution $f:X\to Z$ of $Z$ with irreducible exceptional components $D_i,i\in I$. Let $a_i$ be the discrepancy coefficient of $D_i$. Set $D_J:=\cap_{i\in J} D_i$ for a subset $J\subset I$. Denote the Hodge-Deligne polynomial of $X$ by $\sum_{i,j} a_{i,j} u^iv^j$ and of $D_J,J\neq \emptyset,$ by $\sum_{i,j} a_{i,j}^J u^iv^j$. Then by developing the alternative formulae 1.5 (5) in power series we can write $b_{3,3}$ as \begin{equation*} \begin{split} b_{3,3}= & \ a_{3,3} - \sum_{i\in I} a_{2,2}^{\{i\}} + \sum_{\substack{J\subset I\\ |J|=2}} a_{1,1}^J - \sum_{\substack{J\subset I\\ |J|=3 }} a_{0,0}^J  \\ &+ \sum_{\substack{i\in I\\ a_i =1}} a_{1,1}^{\{i\}} -\sum_{\substack{i\in I\\ a_i =1}} a_{0,0}^{\{i\}} - \sum_{\substack{\{i,j\} \subset I \\ a_i =1 \text{ or } a_j=1}} \delta^{\{i,j\}} a_{0,0}^{\{i,j\}} + \sum_{\substack{i\in I\\ a_i =2}} a_{0,0}^{\{i\}}, \end{split} \end{equation*} where $\delta^{\{i,j\}}\in \{1,2\}$ is the number of components in $\{i,j\}$ with discrepancy 1. The alternating sum on the first line is always nonnegative (this can be shown by the methods of \cite{SchepersVeys}). The problem comes from the term
\[ - \sum_{\substack{\{i,j\} \subset I \\ a_i =1 \text{ or } a_j=1}} \delta^{\{i,j\}} a_{0,0}^{\{i,j\}},\]
which is $-20$ in our example.

\footnotesize{

$\phantom{some place}$

\noindent Jan Schepers\\ Universiteit Leiden\\ Mathematisch Instituut\\ Niels Bohrweg 1 \\ 2333 CA Leiden\\ The Netherlands\\ \emph{E-mail}: jschepers@math.leidenuniv.nl\\

\noindent Willem Veys \\ Katholieke Universiteit
Leuven\\ Departement Wiskunde\\ Celestijnenlaan 200B\\ 3001 Leuven\\
Belgium\\ \emph{E-mail}: wim.veys@wis.kuleuven.be

\end{document}